\documentclass{birkjour}
\usepackage{amsfonts} 
\usepackage{amssymb}
\usepackage{soul} 
 \newtheorem{thm}{Theorem}[section]

 \theoremstyle{definition}
 \newtheorem{defn}[thm]{Definition}
 \theoremstyle{remark}

 \numberwithin{equation}{section}
\begin{document}
%-------------------------------------------------------------------------
% editorial commands: format for journal of geometry 
\title[Ram\rq{}s Theorem]
 {Ram\rq{}s theorem for Trisection}
%----------Author 1----------------------------
\author[R Bhat]{RAMACHANDRA BHAT}
\address{Talent Development Centre, \\ Indian Institute of Science, Challakere \\
Chitradurga - 577536, Karnataka, INDIA. \\
Alumni of: Rameshwara Vidyodaya High School, Itgi, \\
Taluk: Siddapur, Dist: Uttara Kannada, Karnataka, India.}
\email{dsrbhat@gmail.com}
%--------------------------------------------------
\thanks{The support from the Indian Institute of Science is gratefully acknowledged}
\date{February 10, 2019}
\dedicatory{To my Teachers} 
%%% ----------------------------------------------------------------------
\begin{abstract}
While solving problems, if direct methods does not provide solution, indirect methods are explored.  Today, we need an indirect method to solve the problem of angle trisection as the direct methods have been proved not to provide solutions.  The unstoppable curiosity of Geometers and the newer advanced tools available with time have led to newer approaches to progress further.  Results of exploration strictly following the Platonian rules of Euclidean geometry that could help in arriving at the solution is presented here. 
\end{abstract}
\maketitle
%\tableofcontents
\section{Introduction}
In ancient times, there were no formal measurement systems in place.  While working on designs, Greeks defined the dimensions by considering an arbitrary length as a unit of length, for that particular design. They derived other dimensions within that design geometrically by the addition, multiplication, subtraction and division operations. In addition, they knew the right angle and the Pythagoras theorem that helped them to get the square roots of the given line segment.  Therefore, in Euclidean Geometry, one does not measure the dimensions but the expected designs were very precise, accurate and simple to reproduce.  The only geometrical tools used in those days were the ruler (unmarked straightedge) and compasses.  

With the advancement of mathematics and allied subjects, today we could actually make the measurements to show the accuracy and precision of the Greeks and explain them through the algebraic and trigonometric formulations.  However, even today, the Three Famous Problems of Geometry \cite{Wiki, Jagadeeshan} have remained unsolved.  One of them being the trisecting an angle by the Euclidean procedure. 

While many Geometers have been exploring to find solutions to these problems, in 1837, Pierre Wantzel proved that finding a solution to the problem of trisection of an angle is impossible \cite{Wiki, Jagadeeshan}.  This only meant, using only ruler and compasses, it is not possible to trisect an angle of any given value, with the \textit{tools and knowledge available at that instant of time}. It should also be noted that the proof of impossibility considers primarily the constructability of the angle of value equivalent to one-third of the given value and not the trisectability of the given angle directly.  

In the literature,  there are various methods described \cite{Yates, Mallik, Brooks, Forum, Wolfram} to trisect an angle using various additional tools (Marked ruler, trisector, Quadratrix of Hippias, hyperbola, tomahawk, linkages, etc.) and hence, are not Euclidean methods.  However, the only exact and simple procedure demonstrated as well as proved algebraically is by Origami \cite{Richeson, Shima, Kung}, i.e., the Japanese paper folding technique.  However, in these procedures, multiple conditions get satisfied simultaneously in one folding operation. This is the main constraint in the Greek procedure to reproduce. 

\section{The problem definition:} 
The two approaches followed by many geometers in their exploration are: \begin{enumerate}\item[] (i) to divide the given angle into three equal parts, and/or \item[] (ii) to explore the constructability of an angle of value equivalent to one-third of the given angle. 
\end{enumerate}   
For example, some select angles such as 0, 45, 72, 90, 108 and 180 degrees were trisectable \cite{Yates2}, by directly constructing a corresponding angle equivalent to their one-third value and not by trisecting the given angle. This means, \begin{enumerate} \item[(i)] the actual value of the given angle was known (or measured) prior to attempting the trisection, and \item[(ii)] its trisected component is constructible by Euclidean geometry.  \end{enumerate} 

Since the measure of the given angle was unknown, Greeks could only explore the trisectability and not constructibility.  Hence, the definition by the Euclidean geometry may be presented as: 
\begin{defn} 
\textit{\textbf{\lq\lq{}Divide the given angle} (of unknown value) \textbf{into three equal angles} (angular parts) \st{or construct an angle equal to one-third of the given angle (of unknown value?}) using only two tools, viz., (i) an unmarked straight edge (ruler) and (ii) a compasses\rq\rq{}.}
\end{defn}
\section{Approaches considered:} 
Varieties of approaches were considered to reach the goal.  While working out various possibilities, special attention was given \textit{not to compromise on the conditions of Greek construction of using only ruler and compasses}.  In addition, the intention here was to trisect the given angle and not to find the constructability of angle of any given value. Whenever an approach seemed feasible, actual measurement and an algebraic proof was considered to confirm the exactness.  The interesting results obtained during the exploration are given here in the form of a theorem.

\begin{thm}[Ram\rq{}s Theorem:] 
For any given arbitrary angle, there is another unique angle that could be trisected easily by Greek construction (a procedure using only ruler and compasses). 
\end{thm}
\paragraph{To prove the theorem,} let us consider the methods from three different approaches, viz.,  i) Equilateral Triangle Method, ii) Central Angle Method and, iii) Similar Triangles Method. 

\section{Proof of the theorem:}

\subsection{\underline{Method I}: Equilateral Triangle Method} 

\begin{figure}[!ht]  
\includegraphics[scale=0.4]{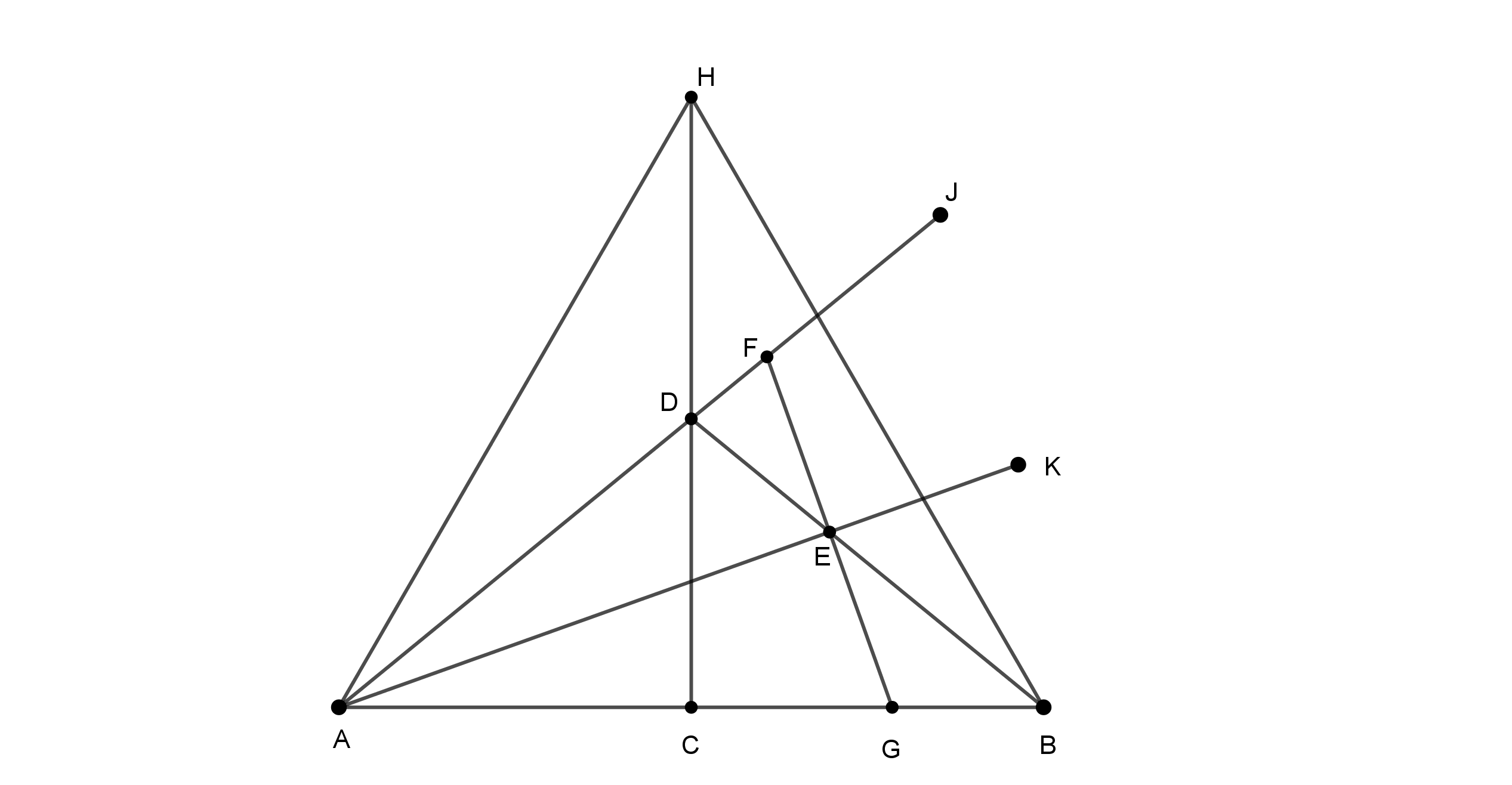}
\caption{Equilateral triangle method.} \label{ETMethod}
\end{figure}

\subsubsection{Procedure for construction (cf Fig. \ref{ETMethod}): }

\begin{enumerate} 
\item[] \textbf{Step 1}: Construct an isosceles triangle with the given angle ($\theta$) at the base. \\ i) Draw a segment AB. Consider this length as 1 (one) unit. \\  ii) Draw a perpendicular bisector for AB. Mark the midpoint of AB as C.  \\ iii) Draw the segments AD and BD such that $\angle$DAC = $\angle$DBC = $\theta$ (given angle). 
\item[] \textbf{Step 2}: Bisect the given angle (at one side). \\  i) Draw the angle bisector for $\angle$DAC. \\ ii) Let it cut the line segment BD at E. Draw the segment AE.  
\item[] \textbf{Step 3}: Draw a perpendicular line to the angle bisector. \\ i) At E, draw a perpendicular line to AE.  \\ ii) This cuts the segment AB at G \& AD (extended) at F. 
\item[] \textbf{Step 4}: Draw an equilateral triangle on the same base. \\ i) Draw an equilateral triangle HAB with AB as the base and H as the third vertex.   
\end{enumerate}

\subsubsection{\textbf{Proof of the theorem:}} 
\begin{enumerate} 
\item[] 
Let $\angle$BAD = $\theta$ (given arbitrary angle). \\ 
$\angle$DBA = $\theta$ (by construction of isosceles triangle). \\
$\angle$BAE = $\angle$EAD = $\frac{\theta}{2}$ (by construction of angle bisector AE). \\
$\angle$FEA = $\angle$AEG = 90$^\circ$ (by construction of perpendicular, AE $\perp$ FG).  \\ 
$\angle$HBA = $\angle$BAH = 60$^\circ$ (by construction of the equilateral triangle). \\
 $\angle$FGA = $\angle$AFG = (90 - $\frac{\theta}{2}$). \\ 
Now, let $\angle$GEB = $\beta$. \\
In triangle BGE, $\angle$GEB + $\angle$EBG = ($\beta$ + $\theta$) = $\angle$EGA = (90 - $\frac{\theta}{2}$). \\ 
Hence, $\angle$GEB = $\beta$ = (90 - $\frac{3\theta}{2}$) or $\theta$ = $\frac{2}{3}$*(90 - $\beta$) = ({60 - $\frac{2}{3}$*$\beta$}) \\
and $\angle$HBE = $\angle$HBA - $\angle$EBA = (60 - $\theta$). \\ 
i.e., $\angle$HBE = (60 - $\theta$) =$\frac{2}{3}$*$\beta$ = $\frac{2}{3}$*$\angle$GEB.\\  Hence, $\angle$HBE is two-third of $\angle$GEB. \\
This implies $\angle$GEB is trisected by the Greek procedure.
\end{enumerate}  

This proves the theorem as we start with an arbitrary angle $\theta$ and end up trisecting a derived angle $\beta$, strictly following the rules of Greek constructions wherein $\theta$ and $\beta$ are related by a unique equation $\beta$ = ({90 - $\frac{3\theta}{2}$}). 

\textbf {Now, by considering the given angle ($\theta$) at $\angle$GEB that is trisected in the process, we need to find the appropriate angle $\angle$BAD to start with.  Hence, the problem of trisection now is to find that unique angle to start with.}

\begin{enumerate}
\item[Note:] i) It may be noted here that when the given angle ($\theta$) happens to be 36$^\circ$, the given angle and the derived angle are of same measure.  Hence, it can be said that the given angle itself is trisected, as a special case. 
\item[] ii) When the starting angle is greater than 60$^\circ$, the value of the derived angle is negative.  This indicates that the trisected component ($\angle$HBE=$\frac{2}{3}\beta$) lies on the outer side of the equilateral triangle (HBA). 
\end{enumerate} 

\subsection {\underline{Method II}: Central Angle Method} 
\begin{figure}[!ht]  
\includegraphics[scale=0.4]{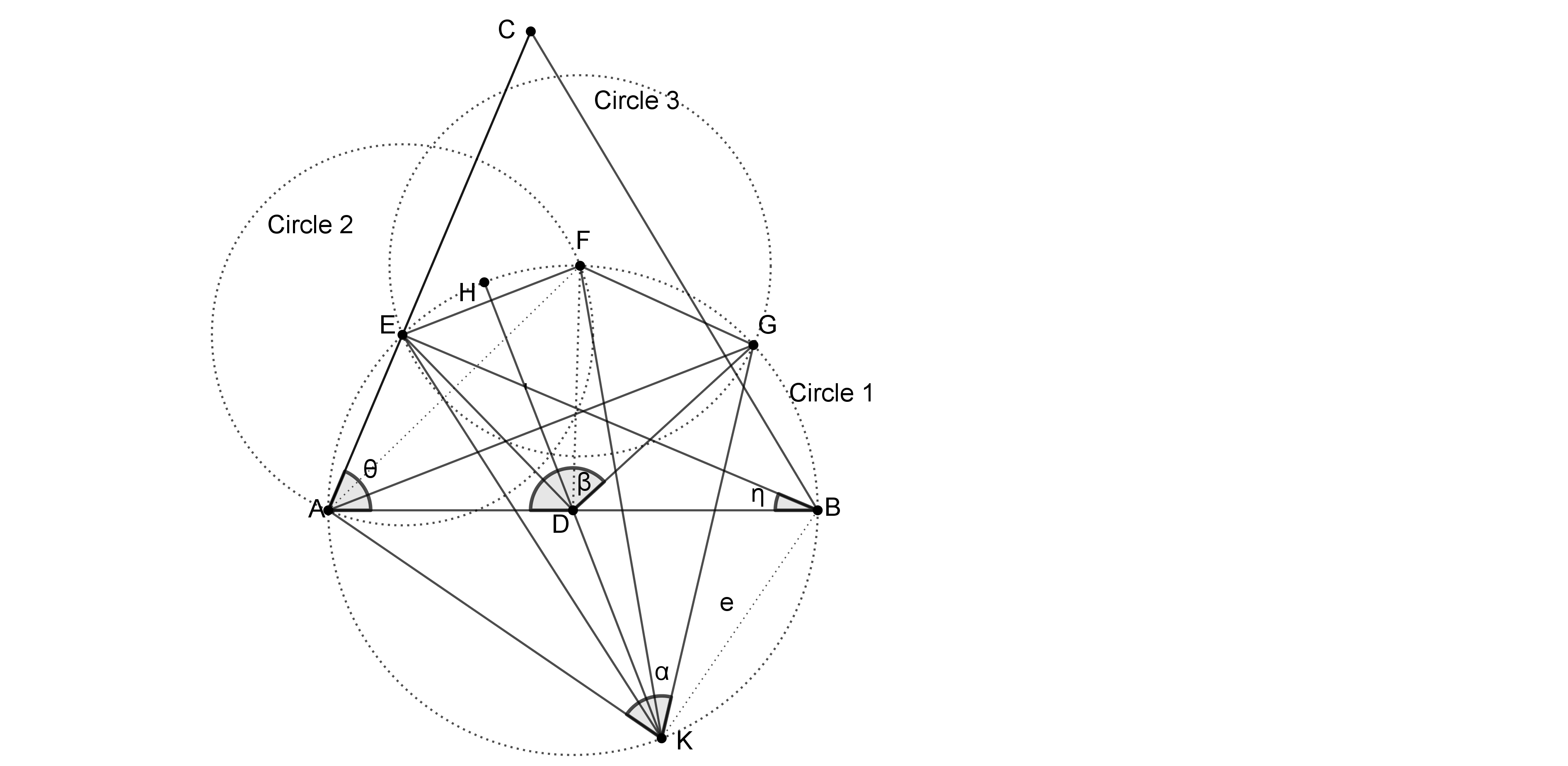} 
\caption{Central Angle Method.} \label{CAMethod} 
\end{figure}

\subsubsection{\textbf {Procedure for construction (cf Fig. \ref{CAMethod}): }} 
\begin{enumerate}
\item[] \textbf{Step 1}: Construct the given angle ($\theta$) with segment AB as the base. \\  i) Draw a segment AB. Consider this length as 1 unit of length.  \\ ii) Draw the given arbitrary angle ($\theta$ = $\angle$CAB). \\ iii) Draw segment AC (=1 unit of length). 
\item[]\textbf{Step 2}: Construct a circle with radius AB as diameter. \\ i) Find the midpoint of AB. Label this midpoint as D. \\  ii) Draw a circle (circle 1) with AD (= DB =$\frac{1}{2}$ unit) as the radius.\\ iii) Label the point of intersection of segment AC and the circle 1 as E.
\item[]\textbf{Step 3}: Construct the equal cords. \\ i) Draw a circle with E as the center and EA as the radius (circle 2). \\ ii) Let F be the other point of intersection of the two circles 1 \& 2. \\ iii) Draw a circle with F as the center and EF as the radius (circle 3). \\  iv) Label the other point of intersection of the two circles 1 \& 3 as G. \\ v) Draw the chords AE, EF, FG, AG, AF, DE, DF, DG and EB. 
\item[]\textbf{Step 4}: Construct an angle at the vertex corresponding to the central angle. \\ i) Draw an angle bisector for $\angle$GDA. \\ ii) Let H and K are the two intersection points of the angle bisector with circle 1. \\ iii) Draw the segments HK, AK, EK, FK and GK. 
\end{enumerate} 
\subsubsection{Proof of the theorem:}
\begin{enumerate}
\item[] 
Let $\angle$EAD = $\theta$ (=given angle) = $\angle$CAB = $\angle$CAD = $\angle$EAB. \\
Let $\angle$GDA = $\beta$; $\angle$EBA = $\eta$; and $\angle$GKA = $\alpha$. \\ 
Now, $\angle$GKA = $\alpha$ = $\frac{1}{2}\angle$GDA = ($\frac{\beta}{2})$  (from the central angle theorem). \\ 
$\angle$AED = $\theta$;   ($\because$ DA = DE = radii; by construction).  \\ 
and $\angle$EBA = $\eta$ = $\frac{1}{2} \angle$EDA = $\frac{1}{2}$(180 - 2*$\theta$) = (90 - $\theta$). \\
Chords AE=EF=FG (by construction) \\
$\therefore$ $\angle$EDA=$\angle$FDE=$\angle$GDF ($\because$ angles of equal arcs). \\ 
But, $\angle$GDA = ($\angle$GDF + $\angle$FDE + $\angle$EDA). \\
$\therefore$  $\angle$EDA = $\angle$FDE = $\angle$GDF = $\frac{1}{3}$ of $\angle$GDA  = $\frac{\beta}{3}$ \\
Similarly, $\    \angle$EKA = $\angle$FKE = $\angle$GKF = $\frac{1}{3}$ of $\angle$GKA = $\frac{\alpha}{3}$  \\
This means $\angle$GDA and $\angle$GKA are trisected by the Greek procedure. \\

This again proves the theorem as we start with an arbitrary angle $\theta$ and end up trisecting a derived angle $\beta$ (and $\alpha$), strictly following the rules of Greek constructions wherein $\theta$ and $\beta$ are related by a unique equation, $\beta$ = 3*(180-2$\theta$) and $\alpha$ = $\frac{\beta}{2}$ = 3*(90-$\theta$) = 3*$\eta$. \\  

Here, if we analyze the diagram (Figure 2) carefully, \\ we could see that $\angle$EAF = $\angle$EKF = ($\frac{\alpha}{3}$) \\ and $\angle$GAF = $\angle$GKF = ($\frac{\alpha}{3}$) (From the central angle theorem). \\
$\angle$BAG = $\frac{1}{2}\angle$BDG = $\frac{1}{2}$(180 - $\beta$) = (90 - $\frac{\beta}{2}$) = (90 - $\alpha$). \\ 
$\angle$BAE=$\angle$BAC=$\theta$=$\angle$BAG+$\angle$GAF+$\angle$EAF=(90 - $\alpha$)+$\frac{\alpha}{3}$+$\frac{\alpha}{3}$=(90 - $\frac{\alpha}{3}$). \\ 

Hence, the relation between the given angle $\theta$ and the derived angle $\beta$ (and $\alpha$) could be re-written as $\theta$=(90 - $\frac{\beta}{6}$)=(90 - $\frac{\alpha}{3}$) . This is the reason why it is said that there is a unique relation between the given angle ($\theta$) and the derived angle ($\beta$ or $\alpha$). 
\end{enumerate}
\begin{enumerate}
\item[Note:] i) Let us consider the difference between the given angle ($\theta$) and the derived angles ($\alpha$) as $\phi$.   Then, $\phi$ = ($\theta$ - $\alpha$) = [(90 - $\frac{\alpha}{3}$) - $\alpha$] = (90 - $\frac{4}{3}*\alpha$).  It should be noted here that when the given angle ($\theta$) happens to be 67.5$^\circ$, i.e., when $\phi$=0, the given angle and the trisected angle are of same measure and hence can be said that the given angle itself is trisected, as a special case. 
\item[] ii) This ‘central angle method’ appears to be a more promising method to find the final solution to the problem.
\end{enumerate} 
\subsection {\underline{Method III}: Similar Triangles Method} 
 
This procedure is primarily an extension of the concept of intersecting secant (angle) theorem \cite{Page} and the diagram used in the algebraic formulation of the problem to derive the trisection equation \cite{Yates3}. 

Let us consider the algebraic formulation of the trisection problem \cite{Yates3} given by the equation x$^{3}$ - 3x - 2a = 0. Here, in the present exploration, the construction (Figure \ref{TEMethod}) is redrawn with an assumed value of a=1 and x=2 and the given arbitrary angle = $3\theta$. \\

\begin{figure}[ht!]  
\includegraphics[scale=0.25]{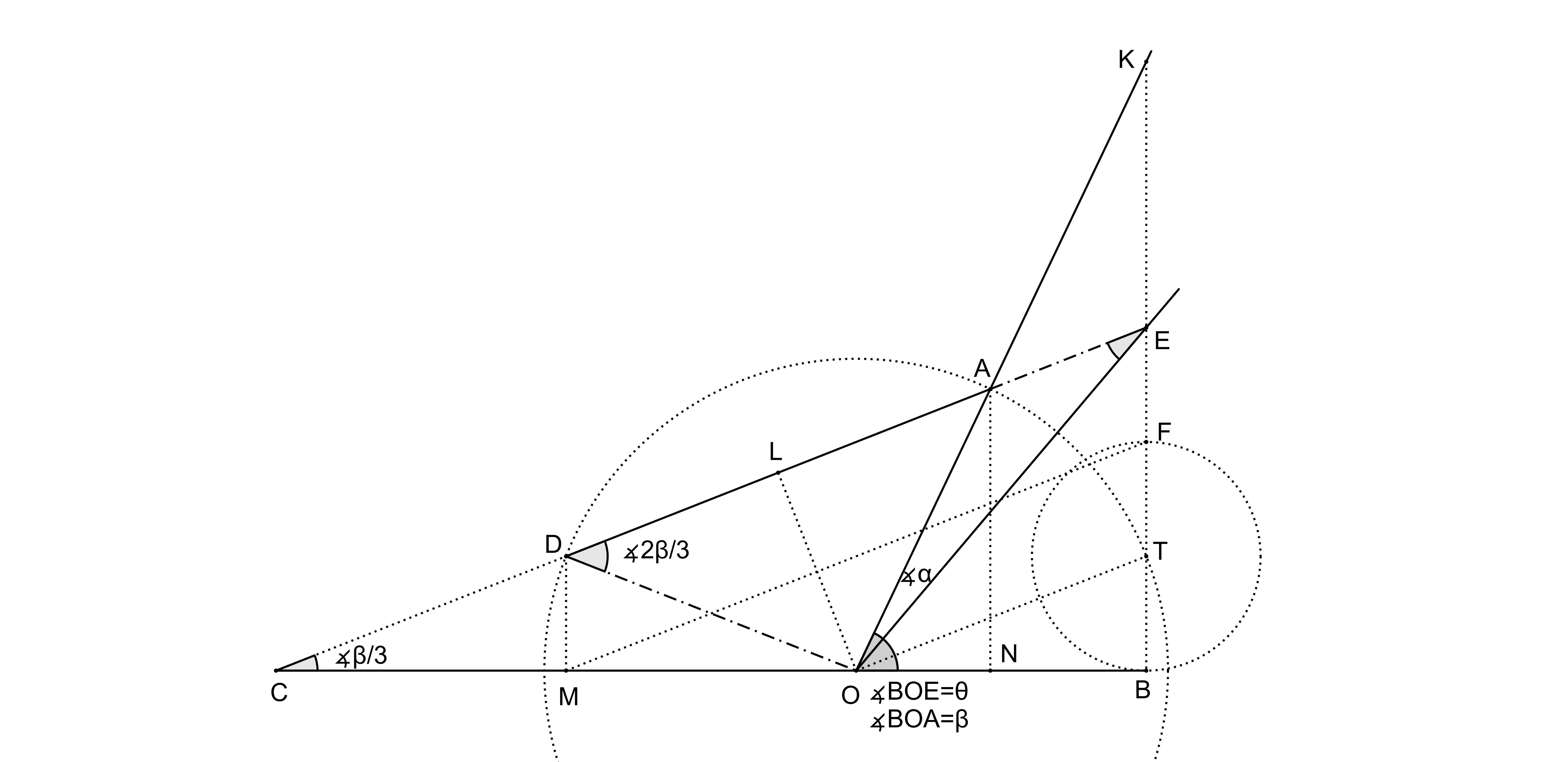} 
\caption{Line Segment Method / Trisection Equation Method.} \label{TEMethod} 
\end{figure}

\begin{enumerate} 
\item[Note:] (i) It should be noted here that these values of a=1 and x=2 are applicable only if 3$\theta$=0$^\circ$ in the algebraic formulation of the trisection equation.  However, in the present study, the author assumed it for any given angle.  \item[] (ii) The given angle is considered as 3$\theta$ (and not $\theta$) for mathematical convenience in \cite{Yates3}. \item[] (iii) In addition, here also the naming of points are kept similar to that in \cite{Yates3}, for easy comparison for the readers.
\end{enumerate} 
\subsubsection{\textbf {Procedure for construction (cf Fig \ref{TEMethod}): }} 
\begin{enumerate}
\item[] \textbf{Step 1}: Construct the given angle ($\theta$).  \\ i) Draw a segment OB. Consider this length (OB) as 1 unit of length.  \\ ii) Draw a perpendicular line to OB at B. \\ iii) Draw a segment OE such that $\angle$BOE is the given arbitrary angle ($\theta$) and, E lies on the perpendicular drawn at B. 
\item[] \textbf{Step 2}: Complete the diagram with the present assumptions (a=1, x=2). \\ i) Extend the segment BO to C such that OC=2*OB (length=2 units). \\ ii) Draw the segments CE \& CO.  Draw a perpendicular bisector to CO. \\ iii) This bisector intersects CE at D.  M is midpoint of CO.  Join MD. \\ iv) Draw a circle with O as center and OD as radius.  \\ v) Mark the intersection of the circle with CE as A \& that of BE as T.  \\ vi) Draw a segment OA and extend it to meet BE (extended) at K. \\ vii) $\angle$BOA=$\beta$. Draw the segment OD [and also CK].  \\ viii) Draw a circle with T as center \& TB as radius. This cut BE at F. Join MF. \\ ix) Draw a perpendicular to CE from O. This cuts CE at L. Join OL. \\ x) Draw a perpendicular to OB from A. It meets OB at N. Draw AN.
\end{enumerate}  

\subsubsection{\textbf {Proof of the theorem}:}
\begin{enumerate} 
\item[] 
The given arbitrary angle $\theta$ = $\angle$BOE and the derived angle $\angle$BOA=$\beta$. \\
From the construction, it is clear that $\angle$MCD = $\frac{\beta}{3}$; $\angle$ODL = $\frac{2\beta}{3}$; $\angle$BOT = $\frac{\beta}{3}$. \\
Here, it may be noted that the proof of the trisection of the angle BOA (=$\beta$) by the above diagram is already available in the literature \cite{Yates}.  Another important point to be noted is that the point T was assumed such that $\angle$BOT = $\frac{\beta}{3}$ (Fig. 2 of ref.\cite{Yates3}), while deriving the trisection equation. However, it is demonstrated here that it is indeed true. The only difference is that, attempt was made to trisect the given angle $\theta$ (=BOE) by the Greek procedure and ended up in trisecting the derived angle $\angle$BOA=$\beta$, with the assumption of OB=1 unit of length. However, the derived angle $\beta$=BOA is uniquely related by the equation $\theta$ = $\beta$+$\alpha$. \\

$\angle$MCD = $\angle$MOD (since CD=OD, OCD is an isosceles triangle). \\
$\therefore$ $\angle$ODL = $\angle$MCD + $\angle$MOD = 2*$\angle$MCD (external angle). \\
$\angle$LAO = $\angle$ODL (since  OD = OA, radii of the same circle). \\ 
$\therefore$  $\angle$LAO=2 *$\angle$MCD.  

Now, $\angle$BOA= $\angle$OCA + $\angle$CAO (external angle).
This is equivalent to $\angle$BOA=$\angle$MCD+$\angle$LAO=$\angle$MCD + 2 *$\angle$MCD = 3 *$\angle$MCD.   \\
Since $\angle$BOA = $\beta$, now $\angle$MCD =$\angle$BOT= $\frac{\beta}{3}$. 
i.e., the derived $\angle$BOA=$\beta$ gets trisected while the starting angle is $\angle$BOE=$\theta$.

Therefore, in this example also, instead of trisecting the given angle ($\theta$), a related angle ($\beta$) gets trisected accurately. This again proves the theorem.

Hence, as an improvement to the current scenario, one could say that any given angle ($\theta$) is trisectable provided one could find a correspondingly related angle from which we should initiate the trisection procedure. 
\end{enumerate} 

\paragraph{Hence, the new problem is not actually trisecting a given angle,} but to find an angle ($\beta$) such that $\theta$ =$\beta$+$\alpha$ that could lead to the trisection of the given angle ($\theta$) easily by the Greek construction.
\section*{Conclusions:} 
It is clear from the three methods given here that the trisection of an angle is possible. However, if the given angle is $\theta$, the angle that is trisected is the derived angle $\beta$ and these two angles ($\theta$ \& $\beta$) are uniquely related. Now, if we want to trisect the given angle ($\theta$) directly, we should place the given angle in the place of the derived angle ($\beta$) as the starting angle and construct back the appropriate angle corresponding to the given angle that could lead to the required trisection by reverse construction. This redefines the problem of trisection as the problem of construction of a unique angle corresponding to the given angle. 

From this study, another spurt of fresh energy among geometers is obvious to explore the problem to conclusion by various approaches. 

While analyzing the various approaches and the methods explored to prove the theorem, it was conjectured that the trisection of the given angle itself could be achieved by the Greek procedure. The same is proposed here as a conjecture. \\ 

\paragraph{Ram’s Conjecture:} 
Any given angle could be trisected by the Greek procedure, 
 i.e., the trisection of an angle (of unknown value) is possible by the Greek procedure.

\subsection*{Acknowledgment}
In the early stage of this study (Mar-2016), while exploring suitable software to verify the accuracy of the ideas, the suggestion of Mr. John Page (mathopenref@gmail.com) to use GeoGebra was of great help and the same is acknowledged here. A special thanks for the creator of the GeoGebra software who made it so user-friendly and available freely on the Internet \cite{Geogebra}.  The support of the Indian Institute of Science while doing this exploration as well as the patience of the family members for the last about four years is greatly acknowledged.

\subsection*{Conflict of Interest}
There is no conflict of interest.

% -----------------
\end{document}